\documentclass[12pt,a4paper]{amsart}
\usepackage{amsfonts,amsmath,amssymb}
\usepackage{hyperref}
\usepackage{latexsym,fullpage,amsfonts,amssymb,amsmath,amscd,graphics,epic}
\usepackage[all]{xy}
\usepackage{amssymb,amsthm,amsxtra}
\usepackage{amscd}
\usepackage{amsthm}
\usepackage{amsfonts}
\usepackage{amssymb}
\usepackage{mathrsfs}
\usepackage{url}
\usepackage{bbm}

\usepackage{enumitem}

\usepackage{}

\newtheorem{theorem}{Theorem}[section]
\newtheorem{lemma}[theorem]{Lemma}
\newtheorem{corollary}[theorem]{Corollary}

\theoremstyle{definition}

\newtheorem{example}[theorem]{Example}

\theoremstyle{remark}

\numberwithin{equation}{section}

\begin{document}

\title{ Quasireductive supergroups}


\author{Vera Serganova}
\address{Dept. of Mathematics, University of California at Berkeley,
Berkeley, CA 94720}
\curraddr{}
\thanks{The work was partially supported by NSF grant
0901554. }

\subjclass[2010]{Primary 17, Secondary 22 }

\date{\today}

\begin{abstract} We call an affine algebraic supergroup quasireductive
if its underlying algebraic group is reductive. We obtain some
results about the structure and representations of reductive supergroups.
\end{abstract}

\maketitle


\section{Introduction}

Reductive algebraic groups have many remarkable applications due to the
fact that all their representations are completely reducible. This
fact, for instance, lies in the foundation of geometric invariant theory.

If one tries to generalize the notion of a reductive group to the category
of algebraic supergroups by imposing the complete reducibility condition,
one immediately discovers that there are very few new objects.
To be precise the only connected  simple reductive supergroups which
are not groups are orthosymplectic supergroups $OSP(1,2n)$.

In these notes we collect preliminary results about algebraic supergroups with
reductive even part. We call such supergroups quasireductive. We hope to
convince the reader that quasireductive supergroups have  manageable
structure and representation theory and there are many interesting open
problems at various level of difficulty in this area.

{\bf Acknowledgments.} The idea to write these notes was inspired by discussions with B. Boe,
C. Boyallian, M. Duflo, R. Fioresi, C. Gruson, P. Heinzner, I. Musson, D. Nakano,
and E. Vishnyakova. Special thanks go to I. Zakharevich for
useful advice and generous help in writing this paper.

\section {Preliminaries}

We work over a field $\mathbb F$  of characteristic not equal to
2. Starting from Section 4 we assume in addition that $\mathbb F$ is
algebraically closed of characteristic zero. A {\it 
superalgebra} is a $\mathbb Z_2$-graded $\mathbb F$-algebra
$A=A_0\oplus A_1$.
All objects and morphisms in the categories of superalgebras and modules
over superalgebras are  $\mathbb Z_2$-graded.
By $p(a)$ we denote the parity of a homogeneous element $a$.

An associative superlagebra $A$ is by definition a $\mathbb Z_2$-graded associative algebra  
and by an $A$-module we always mean a $\mathbb Z_2$-graded
$A$-module. All modules are left unless stated otherwise. The
definitions 
of suprecommutative algebra, tensor product, derivation e.t.c. are
modified from the usual version according to the following supersign rule

{\it all identities are written for homogeneous elements only, and then
extended to all elements by linearity;
whenever in a formula the order of two entries $a$ and $b$ is
switched, the sign $(-1)^{p(a)p(b)}$ appears.}  

A {\it Lie superalgebra} is a vector superspace $\mathfrak g=\mathfrak g_0\oplus\mathfrak g_1$  with an
even( grading preserving) linear map $[\cdot,\cdot]:\mathfrak g\otimes\mathfrak g\to\mathfrak g$, satisfying the following conditions:
      \begin{itemize}
        \item[(i)] $[a,b]=-(-1)^{p(a)p(b)}[b,a]$;
        \item[(ii)] $[a,[b,c]] = [[a,b],c] + (-1)^{p(a)p(b)}[b,[a,c]]$.
      \end{itemize}

By $U(\mathfrak g)$ we denote the universal enveloping algebra of
$\mathfrak g$. For super-analogue of PBW theorem and other
preliminaries on Lie superalgebras we refer the reader to \cite{Kadv}.
By $Z(\mathfrak g)$ we denote the center of $\mathfrak g$ and by
$Z_{\mathfrak g}(\mathfrak h)$ the centralizer of $\mathfrak h$ in $\mathfrak g$. 

\section {Affine algebraic supergroups and Lie superalgebras}

Let $\mathbb F$ be a field and $R=R_0\oplus R_1$ be a commutative Hopf superalgebra
over 
$\mathbb F$ with coproduct $\Delta$ and multiplication $m$. Denote by $I$ the ideal generated by the odd
part $R_1$. Since $I$ is
a Hopf ideal,  
$R/I$ is a commutative Hopf algebra. If $R$ is Noetherian we call $R$
an 
{\it affine supergroup} (sometimes for short we use the term supergroup).
Recall that $R/I$ is the ring $\mathbb F[G_0]$ of regular 
functions on some affine algebraic group $G_0$.
Similarly, one can think about $R$ as the ring of regular functions on
a group object $G$ in the category of affine supervarieties.
Sometimes we will use the notation $R=\mathbb F[G]$.

A Hopf superalgebra $R$ defines a functor $G$ from the category of
unital supercommutative $\mathbb F$-algebras to the category of groups in the
following way
$$G(S)=\{g\in\operatorname{Hom}(R,S)\}.$$
The multiplication on $G(S)$ is defined by
$$gh(r)=g\otimes h(\Delta r)$$
for all $g,h\in G(S)$, $r\in R$, the identity element $e$ coincides
with the counit
$\epsilon$ and the inverse element is given by the formula
$g^{-1}(r)=g(s(r))$, where $s:R\to R$ is the antipode.

We call a supergroup $G$ {\it connected} if $G_0$ is connected and
{\it quasireductive} if $G_0$ is reductive. 

As in the usual case one can define a Lie superalgebra 
$\mathfrak g=\operatorname{Lie}(G)$ as the tangent space to the
supergroup $G$ at the
identity. In the language of Hopf algebras 
$\mathfrak g$ is the space of all $x\in R^*$ satisfying the condition 
$$x(rs)=x(r)\epsilon(s)+(-1)^{p(x)p(r)}\epsilon(r)x(s),$$ 
for all $r,s\in R$. To define a Lie bracket on $\mathfrak g$ for each
$x\in\mathfrak g, r\in R$ set
$$L_x(r)=(\operatorname{id}\otimes x)\circ\Delta(r).$$
Then $L_x:R\to R$  is a derivation of $R$ satisfying the condition
\begin{equation}\label{leftinv}
(\operatorname{id}\otimes L_x)\circ\Delta=\Delta\circ L_x.
\end{equation}
Moreover, every derivation satisfying (\ref{leftinv}) equals $L_x$ for
some $x\in \mathfrak g$. Hence $\mathfrak g$ is closed under
commutator. Thus, $\mathfrak g$ is a Lie superalgebra.
In geometric terms $L_x$ is a left-invariant vector field on $G$.

A left (right) $G$-module $M$ is by definition a right (left)
$R$-comodule. For instance, for
a left $G$-module $M$ the corresponding comodule map
$\alpha:M\to M\otimes R$
satisfies the identity $(\alpha\otimes\operatorname{id})\circ\alpha=(\operatorname{id}\otimes \Delta)\circ\alpha$. 
If $M$ is a left $G$-module, then $M$ has the canonical structure of
a left $\mathfrak g$-module defined by
\begin{equation}\label{module}
xm:=(\operatorname{id}\otimes x)\circ\alpha(m).
\end{equation}

Following \cite{K} and \cite{Kos} we define a {\it  Harish-Chandra pair}
$(\mathfrak g, G_0)$ as the following data 

(i) a Lie superalgebra $\mathfrak g$;

(ii) an algebraic group $G_0$ such that $\operatorname{Lie}(G_0)=\mathfrak g_0$;

(iii) a $G_0$-module structure on $\mathfrak g$ such that the corresponding 
$\mathfrak g_0$-module is adjoint.

Given a Harish-Chandra pair $(\mathfrak g, G_0)$, one can construct a Hopf superalgebra $R=\mathbb F[G]$
such that $\operatorname{Lie}(G)=\mathfrak g$ and $R/I=\mathbb F[G_0]$.
Set $$R:=\operatorname{Hom}_{U(\mathfrak g_0)}(U(\mathfrak g),\mathbb F[G_0]),$$
where the left $U(\mathfrak g_0)$ action on $\mathbb F[G_0]$ is induced by the operators
$L_x$.

Define a multiplication map $m:R\otimes R\to R$ by
$$m(f_1,f_2)(X):=m_0((f_1\otimes f_2)(\Delta_U(X))),$$
where $m_0$ is the multiplication in $\mathbb F[G_0]$ and $\Delta_U$ is the comultiplication
in $U(\mathfrak g)$. It is easy to see that $R$ is a commutative superalgebra isomorphic to
$S(\mathfrak g_1^*)\otimes \mathbb F[G_0]$ (\cite{K}).

Let $\alpha:U(\mathfrak g)\to U(\mathfrak g)\otimes \mathbb F[G_0]$ be the comodule map
such that
$$\alpha(X)_g=\operatorname{Ad}(g^{-1})(X)$$ 
for any $X\in U(\mathfrak g)$ and $g\in G_0$. Let
\begin{equation}\label{adj}
\alpha(X)=\sum_iX_i\otimes r_i.
\end{equation}
Define $\Delta: R\to R\otimes R$ by the formula
\begin{equation}\label{delta}
\Delta f(X,Y):=\sum_i(\operatorname{id}\otimes m_0)\circ(\Delta_0\otimes \operatorname{id})(f(X_iY)\otimes r_i),
\end{equation}
where $\Delta_0$ is the comultiplication in $\mathbb F[G_0]$.
If $g,h\in G_0$, then we have
$$\Delta f(X,Y)_{g,h}=f(\operatorname{Ad}(h^{-1})(X)Y)_{gh}.$$
We define the counit map $\epsilon:R\to\mathbb F$ by
$$\epsilon f =\epsilon_0\circ f(1),$$
where $\epsilon_0$ is the counit in $\mathbb{F}[G_0]$. Finally, define the antipode $s:R\to R$ by
$$sf(X)=s_0\circ m_0\circ (s_0(r_i)\otimes f(s_U(X_i))),$$
where $s_0$ and $s_U$ are the antipodes in  $\mathbb F[G_0]$ and $U(\mathfrak g)$
respectively. If $g\in G_0$, then
$$sf(X)_g=f(\operatorname{Ad}(g)s_U(X))_{g^{-1}}.$$
It is rather tedious job to check that all above operations on $R$ are well defined and that
$R$ is indeed a Hopf algebra satisfying our requirements. We send the reader to the papers
\cite{K},\cite{jap} and \cite{V} for different parts of this
checking. However, we should mention that our 
formulas for $\Delta$ and $S$ are slightly different from those given in \cite{jap} and \cite{V} due to the fact that we use 
the structure of a left $\mathfrak g$-module on $R$. 

The following theorem was proven in \cite{V} for the category of complex analytic supergroups
but can be immediately generalized to the category of algebraic supergroups.

\begin{theorem}\label{exist} 
The category of algebraic Harish-Chandra pairs is equivalent to the
category of algebraic supergroups.
\end{theorem}

By  a $(\mathfrak g,G_0)$-module we understand a $\mathfrak g$-module
$M$ with a compatible $G_0$-module structure.
Let $M$ be a $(\mathfrak g,G_0)$-module and $\rho_0:M\to M\otimes \mathbb F[G_0]$
denote the corresponding $F[G_0]$-comodule map. Define a structure of a $G$-module on $M$ by the following
comodule map $\rho:M\to M\otimes R$ (using the notations of (\ref{adj}))
$$\rho(v)(X)=\sum_i (\operatorname{id}\circ m_0)(\rho_0(X_i v)\otimes r_i)$$  
for any $X\in U(\mathfrak g), v\in M$.
Thus, we obtain a functor from the category of  $(\mathfrak g,G_0)$-modules
to the category of $G$-modules. It is not difficult to see that the inverse functor is given by
(\ref{module}) and therefore there is an equivalence between the categories of  $(\mathfrak g,G_0)$-modules
and $G$-modules (see \cite{Kos}).

\section{Simple Lie superalgebras with reductive even part}

Simple finite-dimensional Lie superalgebras over algebraically closed
$\mathbb F$ of characteristic zero were classified
in \cite{Kadv}.
In this section we review the classification in the case when $\mathfrak g_0$ is reductive.

\subsection{ Basic classical superalgebras}

Let $V$ be a vector superspace of dimension $(m|n)$. The {\it general linear superalgebra}
$\mathfrak {gl}(m,n)$ is by definition the superspace
$\operatorname{End}_\mathbb F(V)$ with the naturally defined bracket.
Let $X$ be an $(m+n)\times(m+n)$-matrix, of the form
      $$
        X =
        \begin{pmatrix}
 
          A&B\\C&D
        \end{pmatrix}.
      $$
Define the \emph{supertrace} of $X$ by
      $$\operatorname{str}(X) = \operatorname{tr}A - \operatorname{tr} D$$
and the {\it special linear Lie superalgebra}  as
$$\mathfrak{sl}(m,n)=\{X\in\mathfrak{gl}(m,n) | \operatorname{str}(X)=0\}.$$
The Lie superalgebra       $\mathfrak{sl}(m,n)$  is simple if and only if $m\neq  n$.
      Otherwise its center coincides with the subalgebra  of scalar matrices, and the quotient superalgebra
      $$
        \mathfrak{psl}(n,n) =  \mathfrak{psl}(n,n)/Z( \mathfrak{psl}(n,n))
      $$
      is simple for $n>1$.

      Fix an even symmetric bilinear form 
      $(\cdot,\cdot)$ on a vector superspace $V$, $\operatorname{dim V} = (m|2n)$ (in the usual
      sense this form is symmetric on $V_0$ and skew-symmetric on $V_1$).
 The {\it orthosymplectic} Lie superalgebra
        $$\mathfrak{osp}(m,2n) = \{x\in\mathfrak{gl}(m,2n)| (xv,w) + (-1)^{p(x)p(v)}(v,xw) = 0\text{ for all homogeneous } v,w \in V\}$$
 is simple if $m,n>0$; the even part of $\mathfrak{osp}(m,2n)$ is isomorphic to
 $\mathfrak{so}(m)\oplus\mathfrak{sp}(2n)$.

The Lie superalgebras $\mathfrak{gl}(m,n)$ and $\mathfrak{osp}(m,2n)$ have an invariant non-degenerate
even symmetric form $(X,Y)=\operatorname{str}(XY)$. It induces the invariant symmetric forms on 
 $\mathfrak{sl}(m,n)$ ($m\neq n$) and  $\mathfrak{psl}(n,n)$.

\subsection{ Two new classical superalgebras}

In the case when $\dim V_0=\dim V_1$ the Lie superalgebra
$\operatorname{End}_\mathbb F(V)$ has
simple subquotients which do not have analogues in the even case.

Assume $V=V_0\oplus V_1$ and $\dim V_0 = \dim V_1 = n$.
      Define
      $$
        \mathfrak{q}(n) =\{
          \begin{pmatrix} A&B\\B&A \end{pmatrix}\in\mathfrak{gl}(n,n)\}
      $$
For any matrix $X\in\mathfrak{q}(n)$ define $\operatorname{otr}X:=\operatorname{tr}B$ and let
$$\mathfrak {sq}(n)=\{X\in \mathfrak{q}(n) | \operatorname{otr}X=0\}.$$ 
It is not hard to see that $\mathfrak{sq}(n)= [\mathfrak{q}(n), \mathfrak{q}(n)]$ if $n>2$. 
By $\mathfrak{pq}(n)$ and $\mathfrak{psq}(n)$ we denote the quotients of
$\mathfrak{q}(n)$ and $\mathfrak{sq}(n)$ respectively by the one-dimensional center.
The superalgebra $\mathfrak{psq}(n)$ is simple if $n>2$.
The odd nondegenerate symmetric form $(X,Y)=\operatorname{otr}(XY)$ is invariant. It induces the
odd invariant form on  $\mathfrak{psq}(n)$.

  Let $(\cdot,\cdot)$ be a nondegenerate odd symmetric form on $V$.
      Then $p(n)$ is the subalgebra of $\mathfrak {gl}(n,n)$ that preserves $(\cdot,\cdot)$.
  The elements of $p(n)$ are block matrices of the following form
      $$\begin{pmatrix} A&B\\C&-A^t \end{pmatrix}$$
  such that $B^t = B, C^t = -C$.
  The commutator 
$$sp(n):=[p(n),p(n)]=p(n)\cap\mathfrak{sl}(n,n)$$
is simple if $n>2$.

\subsection{ Exceptional superalgebras \texorpdfstring{$D(2,1;a)$, $G(1,2)$ and $F(1,3)$}{D(2,1; a), G(1,2) and F(1,3)}}

Let $V_1$ denote the standard $\mathfrak{sl}(2)$-module. 
and $\omega$ denote an $\mathfrak{sl}(2)$-invariant skew-symmetric form on $V_1$. 

The even part of $D(2,1;a)$ is isomorphic to
 $\mathfrak{sl}(2)\oplus\mathfrak{sl}(2)\oplus\mathfrak{sl}(2)$, and the odd part is the exterior tensor
 product $V_1\otimes V_1\otimes V_1$ of three copies of the standard 2-dimensional
 $\mathfrak{sl}(2)$-module $V_1$.    Recall that $S^2 V_1$ is isomorphic to the adjoint
      representation of $\mathfrak{sl}(2)$, denote this isomorphism by
        $\rho: S^2 V_1 \to \mathfrak{sl}(2)$. 
Define a Lie bracket  $D(2,1;a)_1\times D(2,1;a)_1\to D(2,1;a)_0$ by
      \begin{multline*}
        [v_1\otimes v_2\otimes v_3, w_1\otimes w_2\otimes w_3] =
        \alpha \omega(v_1, w_1) \omega (v_2, w_2)  \rho(v_3, w_3)\oplus \\ \oplus
        \beta  \omega(v_1, w_1)  \rho( v_2, w_2)  \omega(v_3, w_3) \oplus
        \gamma \rho( v_1,w_1) \omega(v_2, w_2) \omega (v_3, w_3)
      \end{multline*}
for some $\alpha,\beta,\gamma\in\mathbb F$ such that $a=\frac{\alpha}{\beta}$.
One can check that the super Jacobi identity holds iff $\alpha + \beta + \gamma = 0$.

      Assume that $\alpha + \beta + \gamma = 0$. 
Since the algebra defined above is isomorphic to the one with the triple $(\alpha,\beta,\gamma)$ replaced by 
the $(c\alpha,c\beta,c\gamma)$ for any non-zero $c\in\mathbb F$, each
triple $(\alpha, \beta, \gamma)$ can be associated with a point in $\mathbb P^1$.
The corresponding Lie superalgebra is simple if
$\alpha,\beta,\gamma\neq 0$. In this way we obtain a one parameter
family of  Lie superalgebras 
$D(2,1;a)$ with $a\in\mathbb P^1$.
$D(2,1;a)$ is simple iff  $a\in\mathbb P^1\setminus\{0,-1,\infty\}$.
Furthermore $D(2,1;a)\simeq D(2,1;b)$ iff $a=-1-b$, $a=\frac{1}{b}$ or  $a=\frac{-1}{b+1}$.

The Lie superalgebra $D(2,1;0)$ has the ideal isomorphic to $\mathfrak{psl}(2,2)$ and the quotient 
$D(2,1;0)/\mathfrak{psl}(2,2)$ is isomorphic to $\mathfrak{sl}(2)$.

The exceptional superalgebras $G(1,2)$ and $F(1,3)$  are particular
cases of the following construction. Let $\mathfrak g_0=\mathfrak{sl}(2)\oplus\mathfrak s$,
where $\mathfrak s$ is a simple Lie algebra and $\mathfrak g_1=V_1\otimes V$, where
$V$ is a simple $\mathfrak s$-module.
Assume that there exists  an $\mathfrak s$-invariant
symmetric form $b(\cdot,\cdot)$ on $V$. Then the adjoint $\mathfrak s$-module is a
submodule in $\Lambda^2(V)\subset V\otimes V^*$, where $V^*$ is
identified with $V$ by means of $b$. Hence there is a homomorphism 
$s:\Lambda^2(V)\to\mathfrak s$ of $\mathfrak s$-modules. Define a Lie bracket 
on $\mathfrak g_1\times\mathfrak g_1\to\mathfrak g_0$ by
$$[v\otimes w,v'\otimes w']= \omega(v,
  v')s(w,w')+b(w,w')\rho(v,v')$$
for any $v,v'\in V_1$, $w,w'\in V$.

To construct $G(1,2)$ let $\mathfrak s$ be of type $G_2$ and $V$ be the 
7-dimensional $G_2$-module. 

To construct $F(1,3)$ let $\mathfrak s\simeq\mathfrak{so}(7)$ and $V$ be the spinor
8-dimensional representation of $\mathfrak{so}(7)$.

The Killing form $(x,y)=\operatorname{str}(\operatorname{ad}_x \operatorname{ad}_y)$ is non-degenerate
on all exceptional superalgebras.

\section{Quasireductive Lie superalgebras}

We call a Lie superalgebra $\mathfrak g$ {\it quasireductive} if
$\mathfrak g_0$ is a reductive Lie algebra and
 $\mathfrak g$ is semisimple as a $\mathfrak g_0$-module. It is clear
 that an ideal and a quotient algebra of a quasireductive Lie
 superalgebra is quasireductive. As follows from Theorem \ref{exist} a Lie superalgebra of a quasireductive
supergroup is quasireductive.
A simple Lie superalgebra $\mathfrak g$ with reductive $\mathfrak g_0$ is quasireductive.
The Lie superalgebras $\mathfrak{gl}(n,n)$, $p(n)$ and $\mathfrak q(n)$
are also quasireductive. Below we give more examples of quasireductive Lie algebras.

\begin{example}\label{double} Let $\mathfrak k$ be a simple Lie algebra. Then 
$\mathfrak k^d:=\mathfrak k\otimes \mathbb F(\theta)$ with odd $\theta$ is a 
quasireductive Lie superalgebra with $\mathfrak k^d_0=\mathfrak k^d_1=\mathfrak k$
such that $\mathfrak k^d_1$ is an abelian ideal.

Let $p(\tau)=1$, $p(z_1)=p(z_2)=0$ and 
$$\hat{\mathfrak k}^d=\mathbb F \tau\oplus \mathfrak k^d\oplus \mathbb F z_1\oplus \mathbb F z_2,$$ 
be a vector superspace with a bracket defined by
$$[a\tau+x\otimes 1+y\otimes \theta+b z_1+cz_2,a' \tau+x'\otimes 1+y'\otimes \theta+b' z_1+c'z_2]:=$$
$$([x,x']+a'y+ay')\otimes 1+([x,y']+[y,x'])\otimes \theta+(y,y')z_1+aa'z_2,$$
where $x,x',y,y'\in\mathfrak k$, $a,a',b,b',c,c'\in\mathbb F$ and $(\cdot,\cdot)$ denotes the Killing form on $\mathfrak k$.
It is easy to check that $\hat{\mathfrak k}^d$ is a quasireductive Lie superalgebra. 

Let $$W(0,1)=\mathbb F \frac{\partial}{\partial \theta}\oplus \mathbb F \theta\frac{\partial}{\partial \theta}$$
be the Lie superalgebra of derivations of $\mathbb F(\theta)$. Then
$\mathfrak k^d$ has the obvious structure
of a $W(0,1)$-module. Denote by $\tilde{\mathfrak k}^d$ the semidirect product of  $W(0,1)$ and  $\mathfrak k^d$. 
It is also a quasireductive Lie superalgebra.
\end{example}

\begin{example}\label{qextension} Let 
$$\hat{\mathfrak q}(n)=\mathfrak q(n)\oplus \mathbb Fz$$ 
be a non-trivial central extension of ${\mathfrak q}(n)$
with Lie bracket defined by the formula
$$[x+az,y+bz]=[x,y]+\operatorname{otr}(x)\operatorname{otr}(y),$$
with $x,y\in\mathfrak q(n)$. Then $\hat{\mathfrak q}(n)$ is quasireductive.
Note that $\hat{\mathfrak q}(2)\simeq\hat{\mathfrak{sl}}(2)^d$. 
\end{example}

\begin{example}\label{Clifford} Let $V$ be an $n$-dimensional vector space.
By $\mathfrak {co}(m,n)$ we denote the Lie superalgebra
with the even part $\mathfrak c(m,n)_0=\mathfrak{so}(m)\oplus S^2(V)$ and  
the odd part  $\mathfrak  {co}(m,n)_1=E\otimes V$ such that $S^2(V)$ is the center, $E$ is the standard 
$\mathfrak{so}(m)$-module and the bracket on the odd part is given by
$$[e_1\otimes v_1,e_2\otimes v_2]=b(e_1,e_2)v_1v_2\in S^2(V)$$
for any $v_1,v_2\in V, e_1,e_2\in E$, where $b$ is an $\mathfrak{so}(m)$-invariant symmetric form on $E$.
It is clear that $\mathfrak  {co}(m,n)$ is quasireductive.

Similarly one can define a Lie superalgebra  $\mathfrak {csp}(2m,n)$
with the even part $\mathfrak {csp}(m,n)_0=\mathfrak{sp}(2m)\oplus \Lambda^2(V)$ and  
the odd part  $\mathfrak  {csp}(2m,n)_1=E\otimes V$ such that $\Lambda^2(V)$ is the center, $E$ is the standard 
$\mathfrak{sp}(2m)$-module and the bracket on the odd part is given by
$$[e_1\otimes v_1,e_2\otimes v_2]=\omega(e_1,e_2)v_1\wedge v_2\in \Lambda^2(V)$$
for any $v_1,v_2\in V, e_1,e_2\in E$,  where $\omega$ is an $\mathfrak{sp}(2m)$-invariant skewsymmetric form on $E$.
\end{example}

\begin{example}\label{Clifford2} Let $X$ and $Y$ be vector spaces of dimension
$p$ and $q$ respectively.
By $\mathfrak a(s,p,q)$ we denote the Lie superalgebra
with the even part $\mathfrak a(s,p,q)_0=\mathfrak{gl}(s)\oplus X\otimes Y$ and  
the odd part  $\mathfrak a(s,p,q)_1=E\otimes X+E^*\otimes Y$ such that $X\otimes Y$ is the center, $E$ is the standard 
$\mathfrak{gl}(s)$-module and the bracket on the odd part is given by
$$[e_1\otimes x_1+f_1\otimes y_1, e_2\otimes x_2+f_2\otimes y_2]=f_2(e_1)x_1\otimes y_2+f_1(e_2)x_2\otimes y_1$$
for any $x_1,x_2\in X, y_1,y_2\in Y, e_1,e_2\in E, f_1,f_2\in E^*$.
It is clear that $\mathfrak a(s,p,q)$ is quasireductive.
\end{example}



\begin{lemma}\label{center} If $\mathfrak g$ is quasireductive, then $Z(\mathfrak g/Z(\mathfrak g))_0=0$.
\end{lemma}
\begin{proof} Let $\mathfrak g'=\mathfrak g/Z(\mathfrak g)$ and
$\pi:\mathfrak g\to\mathfrak g'$ be the natural projection.
Let $y\in Z(\mathfrak g')_0$. Choose $x\in Z(\mathfrak g_0)$ such that
$\pi(x)=y$. Then
$\operatorname{ad}_x$ is semisimple. On the other hand,
$\operatorname{ad}_x(\mathfrak g)\subset Z(\mathfrak g).$ Hence $\operatorname{ad}_x(\mathfrak g)=0$.
Thus, $x\in Z(\mathfrak g)$ and $y=0$.
\end{proof}

\begin{lemma}\label{ideal} Let $\mathfrak g$ be quasireductive 
and $\mathfrak l$ be a minimal non-zero ideal in
$\mathfrak g$. Assume that $Z(\mathfrak g)_0=0$.
Then there are the following three possibilities for  $\mathfrak l$

(1)  $\mathfrak l$ is a simple Lie superalgebra;

(2)  $\mathfrak l_1= \mathfrak l$ and $\mathfrak l$ is abelian

(3)  $\mathfrak l=\mathfrak k^d$ for some simple Lie algebra $\mathfrak k$.
\end{lemma}
\begin{proof} Since $\mathfrak g$ is a semisimple $\mathfrak g_0$-module, 
we have a decomposition $\mathfrak g=\mathfrak l\oplus\mathfrak d$ for some
$\mathfrak g_0$-invariant $\mathfrak d$. Then $[\mathfrak l_0,\mathfrak d]\subset \mathfrak d\cap\mathfrak l=0$.

Let $\mathfrak i$ be some ideal of $\mathfrak l$. Assume first
$\mathfrak i_0\neq 0$. It is easy to see that $\mathfrak j=\mathfrak
i_0+[\mathfrak i_0,\mathfrak l]\subset {\mathfrak i}$ is an ideal.
Then we have $[\mathfrak d,\mathfrak i_0]=0$ and 
$$[\mathfrak d,[\mathfrak i_0,\mathfrak l_1]]\subset [\mathfrak
i_0,[\mathfrak d,\mathfrak l_1]]\subset [\mathfrak i_0,\mathfrak
l]\subset \mathfrak j.$$
Therefore $\mathfrak j$ is an ideal in $\mathfrak g$. By minimality of
$\mathfrak l$, $\mathfrak j=\mathfrak l$.

By above every proper ideal of $\mathfrak l$ has trivial even part.
If $\mathfrak l_0=0$, then $\mathfrak l=\mathfrak l_1$ is abelian
as in case(2).

Assume now that $\mathfrak l_0\neq 0$ and $\mathfrak l$ is not simple.   
If  $\mathfrak i\subset \mathfrak l_1$ is an
ideal, then  $[\mathfrak i, \mathfrak i]=0$. If  $\mathfrak s$ is any
ideal of  $\mathfrak l_0$, then  $\mathfrak s\oplus \mathfrak i$ is
also an ideal. Hence $\mathfrak l_0$ does not have non-trivial
proper ideals, $\mathfrak l_1$ is a simple $\mathfrak l_0$-module  and
$[\mathfrak l_1,\mathfrak l_1]=0$. By minimality of $\mathfrak l$, there 
exists $d\in \mathfrak d_1$ such that $[d,\mathfrak l_1]\neq 0$. 
Then the map $d:\mathfrak l_1\to \mathfrak l_0$ is an
isomorphism of $\mathfrak l_0$-modules. Hence $\mathfrak l$ satisfies
the condition (3).
\end{proof}

\section{Derivations and the canonical filtration 
of quasireductive superalgebras}

Let $\operatorname{Der}(\mathfrak g)$ denote the Lie superalgebra of
derivations of $\mathfrak g$ and let $D(\mathfrak g):=\operatorname{Der}(\mathfrak g)/\operatorname{ad}(\mathfrak g)$.

\begin{lemma}\label{der} Let $\mathfrak g$ be quasireductive. Then
$$\operatorname{Der}(\mathfrak g)=\operatorname{Der}^{\mathfrak g_0}(\mathfrak g)+\operatorname{ad}\mathfrak g$$
where $\operatorname{Der}^{\mathfrak g_0}(\mathfrak g)$ denote the
annihilator of ${\mathfrak g_0}$ in $\operatorname{Der}(\mathfrak g)$,
and 
$$D(\mathfrak g)\subset\operatorname{Hom}_{\mathfrak
g_0}(\mathfrak g_1,\mathfrak g)/\operatorname{ad}(Z_\mathfrak g(\mathfrak g_0)).$$
\end{lemma}
\begin{proof} We have a decomposition of $\mathfrak g_0$-modules
$$\operatorname{Der}(\mathfrak g)=\mathfrak d\oplus\operatorname{ad}\mathfrak g.$$
Obviously, 
$\mathfrak d\subset \operatorname{Der}^{\mathfrak g_0}(\mathfrak g)$. That implies the first statement.
The second statement follows from the first immediately.
\end{proof}

\begin{corollary}\label{der1} Let  $\mathfrak g$ be quasireductive.

(1) If $\operatorname{End}_{\mathfrak g_0}(\mathfrak g_1)=Z(\mathfrak g_0)$, then
$D_0(\mathfrak g)=0$.

(2) If $\operatorname{Hom}_{\mathfrak g_0}(\mathfrak g_1,\mathfrak
g_0)=0$, then $D_1(\mathfrak g)=0$.
\end{corollary}

\begin{corollary}\label{ders2}  Let  $\mathfrak l$ be isomorphic to
$\mathfrak{sl}(m,n)$ ($m\neq n$), $\mathfrak{osp}(m,2n)$,
$D(2,1;a)$ ($a\neq 0,-1$), $G(1,2)$ or $F(1,3)$. Then 
 $\operatorname{Der}(\mathfrak l)=\mathfrak l$.
\end{corollary}

\begin{lemma}\label{ders3}  Let $n>2$.  Then  
$\operatorname{Der}(\mathfrak {psl}(n,n))\simeq \mathfrak {pgl}(n,n)$ 
and $\operatorname{Der}(sp(n))\simeq p(n)$. 
\end{lemma}
\begin{proof} Let $\mathfrak l$ be isomorphic to 
$\mathfrak {psl}(n,n)$ or $sp(n)$. 
Note that  $\mathfrak l$ is a type I  superalgebra,
i.e. there exists a $\mathbb Z$-grading 
$$\mathfrak l=\mathfrak l^{-1}\oplus\mathfrak l^{0}\oplus\mathfrak l^{1},$$
with $\mathfrak l^{0}=\mathfrak l_{0}$. Therefore  $\mathfrak l$ has an even derivation $d$ which multiplies
an element by its degree. By direct inspection 
$\operatorname{Hom}_{\mathfrak l_0}(\mathfrak l_1,\mathfrak l_0)=0$,
and the statement follows from Corollary \ref{der1}.
\end{proof}

\begin{lemma}\label{derq}  If $n>2$, then 
$\operatorname{Der}(\mathfrak {psq}(n))=\mathfrak {pq}(n)$ and it is a
semi-direct product of $D(\mathfrak {psq}(n))=\mathbb F^{0|1}$ and $\mathfrak {psq}(n)$.
\end{lemma}
\begin{proof} Let  $\mathfrak l=\mathfrak {psq}(n)$. One can easily
see that 
$\operatorname{Hom}_{\mathfrak l_0}(\mathfrak l_1,\mathfrak l)=\mathbb F^{1|1}$ 
and that a non-zero $d\in\operatorname{End}_{\mathfrak l_0}(\mathfrak l_1)$ does not induce a derivation. So 
$D(\mathfrak l)\subset \mathbb F^{0|1},$
and the lemma follows by Lemma \ref{der}.
\end{proof}
\begin{lemma}\label{der3} Let $\mathfrak l=\mathfrak k^d$ for some simple Lie algebra $\mathfrak k$. Then
$$\operatorname{Der}(\mathfrak l)=\tilde{\mathfrak k}^d$$
and $D(\mathfrak l)\simeq W(0|1)$. 
\end{lemma}
\begin{proof} As in the proof of the previous Lemma we have 
$\operatorname{Hom}_{\mathfrak l_0}(\mathfrak l_1,\mathfrak l)=\mathbb F^{1|1}.$
But in this case the isomorphism $d_1:\mathfrak l_1\to\mathfrak l_0$
induces an odd derivation of $\mathfrak l$ and the map $d_0$ which is
zero on $\mathfrak l_0$ and the identity on $\mathfrak l_1$ is an even
derivation. Now the statement follows by Lemma \ref{der}.
\end{proof}
\begin{lemma}\label{dersl(2)} 
$\operatorname{Der}(\mathfrak {psl}(2,2))=D(2,1;0)$ and $D(\mathfrak {psl}(2,2))\simeq\mathfrak{sl}(2)$.
\end{lemma}
\begin{proof} Let $\mathfrak l=\mathfrak {psl}(2,2)$. 
Obviously, $D(2,1;0)$ is a subalgebra in
$\operatorname{Der}(\mathfrak l)$. By Corollary \ref{der1}(2)
$\operatorname{Der}_1(\mathfrak l)=\mathfrak l_1.$
Furthermore, $\mathfrak l_1$ is a direct sum of two
isomorphic simple $\mathfrak l_0$-modules. Therefore 
$D(\mathfrak l)$ is a Lie subalgebra of 
$\operatorname{End}_{\mathfrak l_0}(\mathfrak l_1)=\mathfrak{gl}(2)$. 
It remains to show that 
$D(\mathfrak l)=\mathfrak{sl}(2)$,
and that follows immediately from the fact that the operator identical
on $\mathfrak g_1$ and zero on $\mathfrak g_0$ is not a derivation.
\end{proof}

\begin{lemma}\label{dersum} Let $$\mathfrak s=\mathfrak m\oplus\mathfrak l^1\oplus\dots\oplus\mathfrak l^k,$$
where $\mathfrak m$ is an odd abelian Lie superalgebra and 
$\mathfrak l^1,\dots,\mathfrak l^k$ satisfy condition (1) or condition
(3) of Lemma \ref{ideal}. 
Then 
$$D(\mathfrak s)=D(\mathfrak m)\oplus D(\mathfrak l^1)\oplus\dots\oplus D(\mathfrak l^k),$$
and $\operatorname{Der}(\mathfrak s)$ is a semidirect-product of
$D(\mathfrak s)$ and $\operatorname{ad}(\mathfrak s)$. 
Moreover, $D_0(\mathfrak s)$ is a reductive Lie algebra and
$D_1(\mathfrak s)$ is an abelian ideal in $D(\mathfrak s)$ trivial as
a module over  $[D_0(\mathfrak s), D_0(\mathfrak s)]$.
\end{lemma}
\begin{proof} We claim that 
$$\operatorname{Der}(\mathfrak s)=\operatorname{Der}(\mathfrak m)\oplus
\operatorname{Der}(\mathfrak
l^1)\oplus\dots\oplus\operatorname{Der}(\mathfrak l^k).$$
Indeed, let $d\in \operatorname{Der}(\mathfrak s)$. By Lemma \ref{der}
we may assume without loss of generality that
$d\in\operatorname{Der}^{\mathfrak s_0}(\mathfrak s)$.
Thus $d:\mathfrak s\to\mathfrak s$ is a homomorphism of $\mathfrak
s_0$-modules. That implies $d(\mathfrak l^i)\subset \mathfrak l^i$, 
$d(\mathfrak m)\subset \mathfrak m$, and the first statement follows.
The statements about a semidirect product and the structure of
$D(\mathfrak s)$ follow from
the similar statements 
for $\mathfrak l^i$ (done above) and for 
$\mathfrak m$ (since $D(\mathfrak m)=\mathfrak{gl}(\mathfrak m)$).
\end{proof}

\begin{theorem}\label{th1} Let $\mathfrak g$ be a quasireductive
Lie superalgebra,  $\mathfrak g':=\mathfrak g/Z(\mathfrak g)$,  $C(\mathfrak g)$ be the sum of all minimal ideals
of $\mathfrak g'$ and $R(\mathfrak g):=\mathfrak g'/C(\mathfrak g)$. 
Then  $C(\mathfrak g)$ 
is a direct sum of ideals satisfying the conditions (1)-(3) of Lemma
\ref{ideal},  $R(\mathfrak g)$ is a semi-direct product of
a reductive even part $\mathfrak r_0$ and the abelian odd ideal with trivial action
of $[\mathfrak r_0,\mathfrak r_0]$. Moreover,  
$\mathfrak g'$ is a semidirect product of 
$R(\mathfrak g)$ and $C(\mathfrak g)$. 
\end{theorem}
\begin{proof}
By Lemma \ref{ideal} and Lemma \ref{center} 
$$C(\mathfrak g)=\mathfrak m\oplus\mathfrak l^1\oplus\dots\oplus\mathfrak l^k,$$
where $\mathfrak m$ is an odd abelian Lie superalgebra and 
$\mathfrak l^1,\dots,\mathfrak l^k$ satisfy condition (1) or condition
(3) of Lemma \ref{ideal}. 

Let $\mathfrak z$ denote the centralizer of
$C(\mathfrak g)$ in $\mathfrak g'$. Then 
$\mathfrak z\cap C(\mathfrak g)=\mathfrak m$. We claim that 
$\mathfrak z=\mathfrak m$. Indeed, let $\mathfrak z=\mathfrak m\oplus
\mathfrak z'$ as a $\mathfrak g'_0$-module. Since $\mathfrak m_0=0$,
$\mathfrak z'$ is an ideal in $\mathfrak g'$. Any minimal ideal inside 
$\mathfrak z'$ lies in $C(\mathfrak g)$. Hence $\mathfrak z'=0$.

Consider the natural homomorphism 
$\varphi:\mathfrak g'\to\operatorname{Der}(C(\mathfrak g))$. Since
$\operatorname{Ker}\varphi=\mathfrak z=\mathfrak m$, the induced
homomorphism $\bar{\varphi}:R(\mathfrak g)\to D(C(\mathfrak g))$ is injective.

By Lemma \ref{dersum} $D(C(\mathfrak g))$ is a
semi-direct product of a reductive Lie algebra and an abelian odd
ideal with the trivial action of the commutator of the even
part. Clearly, the same is true for its quasireductive subalgebra $R(\mathfrak g)$.

Finally, to define an embedding $R(\mathfrak g)\to\mathfrak g'$ let
$\mathfrak r$ be a $\mathfrak g'_0$-submodule such that 
$\mathfrak g'=\mathfrak r\oplus C(\mathfrak g)$. We claim that 
$\mathfrak r$ is a subalgebra. Indeed, it suffices to see that 
$[\mathfrak r_1,\mathfrak r_1]\subset\mathfrak r_0$. But 
$\mathfrak r_1\subset Z_{\mathfrak g'_1}(\mathfrak g'_0)$, the latter
is abelian by Lemma \ref{dersum}. Hence $[\mathfrak r_1,\mathfrak r_1]=0$. 
Since obviously $\mathfrak r\simeq R(\mathfrak g)$, the proof of the theorem is complete.
\end{proof}

By the last theorem any quasireductive Lie superalgebra $\mathfrak g$
has a {\it canonical filtration} by ideals with adjoint
quotients $Z(\mathfrak g)$, $C(\mathfrak g)$ and   $R(\mathfrak g)$.

\begin{corollary}\label{filt} If $\mathfrak g$ is quasireductive, then 
the Loewy length of the adjoint module is not greater than 3.
\end{corollary}
\begin{proof} Use notations of Theorem \ref{th1}. One can reorder
direct summands of $C(\mathfrak g)$ so that
$\operatorname{soc}(\mathfrak g)=Z(\mathfrak g)\oplus\mathfrak l^1\oplus\dots\oplus \mathfrak l^s.$
Denote by $\mathfrak f$ the quotient 
$\mathfrak g/\operatorname{soc}(\mathfrak g).$
By Theorem \ref{th1} $\mathfrak f$ is a semi-direct product of the ideal
$\mathfrak i=\mathfrak m\oplus\mathfrak l^{s+1}\oplus\dots\oplus\mathfrak l^k$ 
and the subalgebra $\mathfrak r\simeq R(\mathfrak g)$. 

We claim that
$[[\mathfrak r_1,\mathfrak r_0],\mathfrak i]=0$. Indeed, let $d\in [\mathfrak r_1,\mathfrak r_0]$. 
There exists $c\in\mathfrak r_0$ such that $[c,d]=d$. By the parity argument,
$[d,\mathfrak m]=0$. Let $[d,\mathfrak l^i]\neq 0$. Then  $[c,\mathfrak l^i]\neq 0$, $\mathfrak l^i\simeq \mathfrak k^d$, and 
as follows from Lemma \ref{der3}, there exists $t\neq 0$ such that 
\begin{equation}\label{eq1}
[c,x]=tx\,\,\text{for any}\,\, x\in\mathfrak l^i_1.
\end{equation}
Let $\pi:\mathfrak g\to \mathfrak g'$ be the canonical projection, then we have
the decomposition of $\mathfrak g_0$-modules
$\pi^{-1}(\mathfrak l^i)=\mathfrak l^i\oplus Z(\mathfrak g)$. 
But (\ref{eq1}) ensures that $[\mathfrak l^i_1,\mathfrak l^i_1]=0$. Therefore 
$\mathfrak l^i\subset \operatorname{soc}(\mathfrak g)$ and $i\leq s$.
Thus, the filtration 
$$\operatorname{soc}(\mathfrak g)\subset \pi^{-1}(C(\mathfrak g)\oplus [\mathfrak r_1,\mathfrak r_0])\subset \mathfrak g$$
has semisimple adjoint quotients and the Loewy length of $\mathfrak g$ is at most 3.
\end{proof}

\section {Central extensions and a general construction of a quasireductive
Lie superalgebra}

Our next step is to study central extensions in the category of
quasireductive Lie algebras.  First we note that any exact sequence
$$0\to\mathbb F^{0|1}\to\hat {\mathfrak g}\to\mathfrak g\to 0$$ 
of quasireductive superalgebras splits since it splits over $\mathfrak g_0$.
Therefore a non-trivial central extension in the category of quasireductive Lie
superalgebras is even. Any such extension is described by an even 2-cocycle
$c\in\Lambda^2(\mathfrak g^*)$ such that $c(\mathfrak g_0,\mathfrak g)=0$. It is not difficult to see that  
$c$ satisfies the cocyle condition $dc=0$ iff
$c\in (\Lambda^2(\mathfrak g^*))^{\mathfrak g_0}$. On the other hand,
$c$ defines a trivial central extension iff it satisfies the coboundary
condition, i.e. there exists $f\in \mathfrak g_0^*$ such that 
$c(x,y)=f([x,y])$ for any $x,y\in\mathfrak g$. We denote by
$H^2_r(\mathfrak g)$ the abelian group of non-trivial central
extensions of $\mathfrak g$ in the category of quasireductive
superalgebras. The above discussion implies that
\begin{equation}\label{ext}
H^2_r(\mathfrak g)=(\Lambda^2(\mathfrak g_1^*))^{\mathfrak g_0}/(\mathfrak g_0^*/[\mathfrak g,\mathfrak g]_0^\perp).
\end{equation}
In particular,
\begin{equation}\label{dimext}
\operatorname{dim} H^2_r(\mathfrak g)=\operatorname{dim}(\Lambda^2(\mathfrak g_1^*))^{\mathfrak g_0}-
\operatorname{dim}(Z(\mathfrak g_0)\cap[\mathfrak g_1,\mathfrak g_1]).
\end{equation}

\begin{lemma}\label{centext}

(a) If $\mathfrak l$ is isomorphic to
$\mathfrak{sl}(m,n)$ ($m\neq n$), $\mathfrak{osp}(m,2n)$,
$D(2,1;a)$ ($a\neq 0,-1$), $sp(n)$ ($n>2,n\neq
4$), $G(1,2)$  or $F(1,3)$, then $H^2_r(\mathfrak l)=0$.

(b) If $\mathfrak l$ is isomorphic to
$\mathfrak{psl}(n,n)$ ($n>2$), $\mathfrak{psq}(n)$ ($n>2$) or
$sp(4)$, then  $H^2_r(\mathfrak l)=\mathbb F$. If $\mathfrak l=\mathfrak{psl}(n,n)$ ($n>2$) the corresponding non-trivial central
extension is isomorphic to $\mathfrak{sl}(n,n)$.
If $\mathfrak l=\mathfrak{psq}(n,n)$ ($n>2$) the corresponding non-trivial central
extension is isomorphic to $\mathfrak{sq}(n,n)$.
We denote by $\hat{sp}(4)$ the central extension of 
$sp(4)$. 

(c) If $\mathfrak l=\mathfrak k^d$ for some simple Lie algebra $\mathfrak k$, then  $H^2_r(\mathfrak l)=\mathbb F$.

(d) If $\mathfrak l=\mathfrak {psl}(2,2)$, then  $H^2_r(\mathfrak l)=\mathbb F^3$. 
By $\hat{\mathfrak{psl}}(2,2)$ we denote the 3-dimensional universal central extension of $\mathfrak {psl}(2,2)$. 
\end{lemma} 
\begin{proof} The proof follows from (\ref{dimext}) by direct inspection.
\end{proof}
\begin{lemma}\label{simpleh2} If $\mathfrak l$ is a simple quasireductive superalgebra, then
$H^2(\mathfrak l)=H_r^2(\mathfrak l)$, i.e. every central extension of $\mathfrak l$ is quasireductive.
\end{lemma}
\begin{proof} Since $\operatorname{dim}Z(\mathfrak l_0)\leq 1$, $H^2(\mathfrak l_0)=0$. That implies 
$H^2(\mathfrak l)_0=H_r^2(\mathfrak l)_0$. On the other hand, every odd cocycle $c\in H^2(\mathfrak l)_1$ is $\mathfrak l_0$-invariant. 
If $\mathfrak l\neq\mathfrak{psq}(n)$, then $(\mathfrak l_0\otimes\mathfrak l_1)^{\mathfrak l_0}=0$ and therefore $H^2(\mathfrak l)_1=0$.
Now assume that $\mathfrak l=\mathfrak{psq}(n)$, then an odd cocycle must be of the form
$c(x,y)=t\operatorname{tr}(xy)$, where $x\in\mathfrak l_1, y\in\mathfrak l_0, t\in\mathbb F$. However, there exists
$x\in\mathfrak l_1$ such that $\operatorname{tr}x^3\neq 0$. Then  $c(x,[x,x])=2t\operatorname{tr}x^3$. But $c(x,[x,x])=0$ by the cocycle condition.  
Hence $c=0$ and we have $H^2(\mathfrak l)_1=H_r^2(\mathfrak l)_1=0$.
\end{proof}
\begin{lemma}\label{extsum} Let $$\mathfrak s=\mathfrak m\oplus\mathfrak l^1\oplus\dots\oplus\mathfrak l^k,$$
where $\mathfrak m$ is an odd abelian Lie superalgebra and 
$\mathfrak l^1,\dots,\mathfrak l^k$ satisfy condition (1) or condition
(3) of Lemma \ref{ideal}. Then 
$$H^2_r(\mathfrak s)=H^2_r(\mathfrak m)\oplus H^2_r(\mathfrak
l^1)\oplus\dots\oplus H^2_r(\mathfrak l^k).$$
\end{lemma}
\begin{proof} Follows from (\ref{ext}) and the fact that 
$$(\mathfrak l^i_1\otimes \mathfrak m)^{\mathfrak s_0}=(\mathfrak l^i_1\otimes \mathfrak l^j_1)^{\mathfrak s_0}=0$$
for all $i\neq j$.
\end{proof}

 We call a Lie superalgebra $\mathfrak g$ {\it reduced} if 
$Z(\mathfrak g)\subset [\mathfrak g,\mathfrak g]$. It is clear that any Lie superalgebra is a direct sum
of an abelian Lie superalgebra and a reduced Lie superalgebra.

\begin{corollary}\label{str} Let $\mathfrak g$ be a quasireductive reduced Lie superalgebra, 
$C(\mathfrak g)=\mathfrak l^1\oplus\dots\oplus\mathfrak l^k\oplus \mathfrak m$ such that $\mathfrak m$ is odd abelian, 
each $\mathfrak l_i$ is either simple or $\mathfrak k^d$ for some simple Lie algebra $\mathfrak k$ and  
$\mathfrak r:=R(\mathfrak g)\subset D(\mathfrak s)$. There exists a subspace 
$$Z\subset H^2_r( \mathfrak g')=\bigoplus_{i=1}^kH^2_r(\mathfrak l_i)^{\mathfrak r_0}\oplus \Lambda^2(\mathfrak m^*\oplus \mathfrak r_1^*)^{\mathfrak r_0},$$ 
such that $\mathfrak g$ is the central extension   
$$0\to Z^*\to\mathfrak g\to\mathfrak g'\to 0,$$
where the cocycle $c:\Lambda ^2(\mathfrak g')\to Z^*$ is 
defined by the formula $c(x,y)(z)=z(x,y)$ for any $x,y\in  \mathfrak g', z\in Z$.
Furthermore the above construction defines a reduced Lie superalgebra iff 
$\operatorname{Ker} c\cap \mathfrak m=0$. 
\end{corollary}

\begin{proof} Corollary follows from Theorem \ref{th1}, Lemma \ref{extsum} and (\ref{ext}).
\end{proof} 

\section{On classification of some quasireductive Lie superalgebras}

It is more or less obvious from Corollary \ref{str} that a complete
classification of quasireductive superalgebras is a wild problem. To list all quasireductive reduced
$\mathfrak g$ with fixed $C(\mathfrak g)$ we have to choose a reductive $R(\mathfrak g)$ in  $D(C(\mathfrak g))$
and then a subspace of  $R(\mathfrak g)_0$-invariants in  $H^2_r( \mathfrak g')$. Such description will have infinite
number of parameters and it does not provide us with better
understanding of the situation. The classification of quasireductive superalgebras 
which have an even non-degenerate invariant symmetric form is obtained in \cite{B}.
Here we discuss some other specific cases when some classification can be done.

Introduce the partial order $\leq$ on the set of isomorphisms classes of quasireductive reduced Lie superalgebras
as follows. We say that $\mathfrak f\leq\mathfrak g$ if $C(\mathfrak f)=C(\mathfrak g)$ and there exist a central subalgebra $Z\subset\mathfrak g$
and an embedding $\mathfrak f\to\mathfrak g/Z$. 
A {\it maximal} quasireductive Lie superalgebra is by definition maximal with respect to this order.
Theorem \ref{th1} and Corollary \ref{str} imply the following lemma.

\begin{lemma}\label{max} Recall the notations of Corollary \ref{str}. 
A  quasireductive Lie
 superalgebra $\mathfrak g$ is maximal if and only if 
 $Z=H^2_r(\mathfrak g')^{R(\mathfrak g)_0}$ and 
$R(\mathfrak g)_0=\operatorname{Ann}_{D(\mathfrak g)_0}(Z)$. 

There exists at most one maximal quasireductive superalgebra with fixed $C(\mathfrak g)$
and fixed $R(\mathfrak g)_0$.
\end{lemma}

We call a simple Lie superalgebra $\mathfrak l$ {\it rigid} if $H^2(\mathfrak l)=D(\mathfrak l)=0$. The following simple quasireductive 
Lie superalgebras are rigid: $\mathfrak{sl}(m,n)$ ($m\neq n$), $\mathfrak{osp}(m,2n)$,
$D(2,1;a)$ ($a\neq 0,-1$), $G(1,2)$, $F(1,3)$ and any simple Lie algebra (see Corollary \ref{ders2} and Lemma \ref{centext}). 

We call a quasireductive Lie superalgebra $\mathfrak g$ {\it elementary} iff $C(\mathfrak g)$ is simple or isomorphic 
to $\mathfrak k^d$ for some simple Lie algebra $\mathfrak k$.
\begin{lemma}\label{elementary} Let $\mathfrak g$ be a maximal elementary superalgebra. Then $\mathfrak g$ is rigid simple
or isomorphic to one of the following list: $\mathfrak{gl}(n,n)$ ($n\geq 2$),  $p(n)$ ($n>2$),  $\hat{sp}(4)$, 
$D(2,1;0)$,  $\hat{\mathfrak{psl}}(2,2)$, $\tilde {\mathfrak{q}}(n)$ ($n>2$),  $\hat{\mathfrak{k}}^d$ or  $\tilde{\mathfrak{k}}^d$ for 
some simple Lie algebra $\mathfrak k$. 
\end{lemma}
\begin{proof} We notice first that if $C(\mathfrak g)$ is rigid simple then obviously $C(\mathfrak g)=\mathfrak g$.
Now we consider case by case all other possibilities for $\mathfrak s=C(\mathfrak g)$.

Observe that if  $H^2_r(\mathfrak s)^{D(\mathfrak s)}=H^2_r(\mathfrak s)$, then there is exactly one maximal superalgebra
with $C(\mathfrak g)=\mathfrak s$ and it is the universal central extension of $\operatorname{Der}(\mathfrak s)$. That takes care of the cases
$\mathfrak s=sp(n)$ ($n>2, n\neq 4$), $\mathfrak s=\mathfrak{psl}(n,n)$ ($n>2$) and  $\mathfrak s=\mathfrak{psq}(n)$ ($n>2$).

Let $\mathfrak s=sp(4)$. It is easy to check  that $p(4)$ does not have central extensions,
hence there are two possibilities:  $\mathfrak g=\operatorname{Der}(sp(4))=p(4)$
or  $\mathfrak g=\hat{sp}(4)$. 

Next  consider the case $\mathfrak s=\mathfrak{psl}(2,2)$. Then $D(\mathfrak s)=\mathfrak{sl}(2)$. We have the following three cases

(1) If $R(\mathfrak g)=\mathfrak{sl}(2)$, then  $H^2_r(\mathfrak s)^{\mathfrak {sl}(2)}=0$ and $\mathfrak g=D(2,1;0)$. 

(2) If $R(\mathfrak g)=\mathbb F$, then $\mathfrak
g'=\mathfrak{pgl}(2,2)$ and  $\mathfrak g$ is its universal central
extension isomorphic to $\mathfrak{gl}(2,2)$.

(3) If $R(\mathfrak g)=0$, then  $\mathfrak g$ is isomorphic to the universal central extension $\hat{\mathfrak{psl}}(2,2)$.   

Finally, let  $\mathfrak s=\mathfrak{k}^d$. Then $D(\mathfrak s)=W(0,1)$. If $R(\mathfrak g)=D(\mathfrak s)$ or $D(\mathfrak s)_0$, then
 $H^2_r(\mathfrak s)^{R(\mathfrak s)}=0$ and by maximality $\mathfrak g=\operatorname{Der}(s)=\tilde{\mathfrak{k}}^d$. If
 $R(\mathfrak g)=D(\mathfrak s)_1$ or $0$, then $H^2_r(\mathfrak s)^{R(\mathfrak s)}=\mathbb F$. Hence $\mathfrak g= \hat{\mathfrak{k}}^d$.
\end{proof}

 We call a quasireductive Lie superalgebra $\mathfrak g$ {\it pseudoabelian}
if  $C(\mathfrak g)$ is abelian.

\begin{lemma}\label{pseudoabelian} Every maximal pseudoabelian 
superalgebra is isomorphic to a direct sum 
\begin{equation}\label{abelian}
\bigoplus\mathfrak{co}(m_i,n_i)\oplus \bigoplus \mathfrak{csp}(2r_j,l_j)\bigoplus \mathfrak a(s_k,p_k,q_k),
\end{equation}
with $p_kq_k\neq 1$.
\end{lemma}

\begin{proof} Let $\mathfrak g$ be as in (\ref{abelian}).
We leave it as an exercise for the reader to check that the condition
of Lemma \ref{max} holds for $\mathfrak g$.

Now we will prove that for any pseudoabelian superalgebra $\mathfrak f$ there exists $\mathfrak g$ as above such that
$\mathfrak f\leq\mathfrak g$. Let $\mathfrak m=C(\mathfrak f)$, $\mathfrak r=R(\mathfrak f)$. Write
$$\mathfrak m= \bigoplus_{i}L_i\otimes V_i \oplus \bigoplus_{j}N_j\otimes W_j\oplus \bigoplus_k (M_k\otimes X_k\oplus M^*_k\otimes Y_k),$$
where $L_i$ are non-isomorphic simple $\mathfrak r$-components with invariant symmetric form,
$N_j$ are non-isomorphic simple $\mathfrak r$-components with invariant skew-symmetric form,
and $M_k$ are non-isomorphic
simple components such that  $M_k$ is not isomorphic to  $M_k^*$.
Let $m_i,n_i,2r_j,l_j,s_k,p_k,q_k$ denote the dimensions of $L_i,V_i,N_j,W_j,M_k,X_k,Y_k$ respectively.
Then $\mathfrak f\leq\mathfrak g$.
\end{proof} 

\section{Representation theory}

\subsection{General properties}
Let $G$ be a quasireductive connected supergroup and $\mathfrak g$ be the Lie superalgebra
of $G$. Let $\operatorname{Rep}(G)$ denote the category of all representations of $G$ and 
$\mathcal F(G)$ denote the category of finite-dimensional representations. By
 $\mathcal F(\mathfrak g)$ we denote the category of finite-dimensional $\mathfrak g$-modules
semi-simple over $\mathfrak g_0$. By Section 3 $\mathcal F(G)$ is the full subcategory of $\mathcal F(\mathfrak g)$. 
We require that morphisms between modules preserve parity. So the space 
$\operatorname{Hom}_{\mathfrak g}(M,N)$ is even by definition. By $\Pi$ we denote the functor that changes the parity.

\begin{lemma} \label{inj}(\cite{BKN}) Every module in  $\mathcal F(\mathfrak g)$ has an injective hull and a projective cover.
\end{lemma}
\begin{proof} The category  $\mathcal F(\mathfrak g_0)$ is semisimple by definition. We have two functors
$\operatorname{Ind}:\mathcal F(\mathfrak g_0)\to \mathcal F(\mathfrak g)$ and
$\operatorname{Coind}:\mathcal F(\mathfrak g_0)\to \mathcal F(\mathfrak g)$,
defined by
$$\operatorname{Ind}(M)=U(\mathfrak g)\otimes_{U(\mathfrak g_0)}M$$
and
$$\operatorname{Coind}(M)=\operatorname{Hom}_{U(\mathfrak g_0)}(U(\mathfrak g),M).$$
The first functor maps projective modules to projective modules, the second maps
injective to injective. If $M\in\mathcal F(\mathfrak g)$ we have the natural projection
$\operatorname{Ind}(M)\to M$ and the natural embedding $M\to \operatorname{Coind}(M)$.
That proves the statement.
\end{proof}

For arbitrary module $M\in\mathcal F(\mathfrak g)$ we denote by $I(M)$ and $P(M)$ the injective hull
and the projective cover of $M$ respectively. If $M$ is simple, then the algebras
$\operatorname{End}_\mathfrak g(I(M))$ and
$\operatorname{End}_\mathfrak g(P(M))$ are local. Hence $I(M)$ is an
indecomposable module with unique simple submodule $M$ and $P(M)$ is
an indecomposable module with unique simple quotient $M$.

\begin{theorem}\label{PW} We have the following decomposition of left $\mathfrak g$-modules
$$\mathbb F[G]=\bigoplus I(L)^{\oplus\operatorname{dim}L_0},$$
where summation is taken over all mutually non-isomorphic simple
$L\in\mathcal F(G)$.
\end{theorem}
\begin{proof} As follows from Section 3, 
$$\mathbb F[G]=\operatorname{Coind}(\mathbb F[G_0]).$$ So 
$\mathbb F[G]$ is injective and  
$$\mathbb F[G]=\bigoplus I(L)^{\oplus m(L)},$$ where 
$m(L)=\operatorname{dim}\operatorname{Hom}_{\mathfrak g}(L,\mathbb F[G])$.

Since the group $G_0$ is reductive
$$\mathbb F[G_0]=\bigoplus_S S^{\oplus\operatorname{dim}S},$$
where summation is taken over the set of isomorphism classes of simple
even $G_0$-modules. Now we calculate $m(L)$ using Frobenius reciprocity 
$$m(L)=\operatorname{dim}\operatorname{Hom}_{\mathfrak
g}(L,\operatorname{Coind}(\mathbb F[G_0]))=
\operatorname{dim}\operatorname{Hom}_{\mathfrak g_0}(L,\mathbb F[G_0])=$$
$$\sum_S\operatorname{dimHom}_{\mathfrak g_0}(L,S)\operatorname{dim}S=\operatorname{dim}L_0.$$
\end{proof}

{\bf Problem.} It is very interesting to understand the structure
of $\mathbb F[G]$ as $G\times G$-module. In case $G=GL(m,n)$ it is
done in \cite{SZ}.

\begin{lemma}\label{twist} The even linear functional $f:\mathfrak g\to\mathbb F$ defined by
$$f(x)=\operatorname{tr}_{\mathfrak{g}_1}(\operatorname{ad}_x)$$
for any $x\in\mathfrak g_0$ defines a character  of $\mathfrak g$. 
\end{lemma}
\begin{proof} All we have to check is that $f([y,z])=0$ for any $y,z\in\mathfrak g_1$.
We have $\operatorname{str}(\operatorname{ad}_{[y,z]})=0$. Since
$$\operatorname{str}(\operatorname{ad}_{x})=\operatorname{tr}_{\mathfrak{g}_0}(\operatorname{ad}_x)-\operatorname{tr}_{\mathfrak{g}_1}(\operatorname{ad}_x)$$
and $\operatorname{tr}_{\mathfrak{g}_0}(\operatorname{ad}_x)=0$ for any $x\in\mathfrak g_0$ (because $\mathfrak g_0$ is reductive),
we obtain the statement.
\end{proof}

Let $T$ be the one dimensional $\mathfrak g$-module defined by 
the character $f$ with the same parity as $\operatorname{dim}\mathfrak g_1$. 
If $\mathfrak g$ is not quasireductive $T$ is just a $\mathfrak g_0$-module. 

\begin{lemma}\label{reciprocity} Let $\mathfrak g$ be any Lie superalgebra and $M$ be a $\mathfrak g_0$-module.
Then
$$\operatorname{Coind}(M\otimes T)\simeq\operatorname{Ind}(M).$$
If $\mathfrak g$ is quasireductive then 
$$\operatorname{Coind}(M)\otimes T\simeq\operatorname{Ind}(M).$$
\end{lemma}
\begin{proof} By PBW theorem we have an isomorphism 
$U(\mathfrak g)=U(\mathfrak g_0)\otimes S(\mathfrak g_1)$. 
Note that $T\simeq  S^d(\mathfrak g_1)$ as a $\mathfrak g_0$-module, where
$d=\operatorname{dim}\mathfrak g_1$. Let $\pi: S(\mathfrak g_1)\to  S^d(\mathfrak g_1)$
denote the natural projection.
There are canonical homomorphisms
$$\varphi\in\operatorname{Hom}_{\mathfrak g_0}(M,\operatorname{Coind}(M\otimes T)),$$
defined by 
$\varphi(m)(X):=m\otimes \pi(X)$ for any $m\in M, X\in  S(\mathfrak g_1)$,
and
$$\psi\in\operatorname{Hom}_{\mathfrak g_0}(\operatorname{Ind}(M),M\otimes T),$$
defined by $\psi (X\otimes m):=m\otimes\pi(X)$.
By Frobenius reciprocity they induce mutually inverse homomorphisms 
$\bar{\varphi}:\operatorname{Ind}(M)\to \operatorname{Coind}(M\otimes T)$
and $\bar{\psi}:\operatorname{Coind}(M\otimes T)\to \operatorname{Ind}(M)$. 

The second statement follows immediately from the first one and Lemma \ref{twist}.
\end{proof}

\begin{lemma}\label{duality} Let $\mathfrak g$ be quasireductive and
$L\in\mathcal F(\mathfrak g)$ be simple,
then $I(L\otimes T)\simeq P(L)$.
\end{lemma}
\begin{proof} Any projective module $P$ in $\mathcal F(\mathfrak g)$
has a direct sum decomposition
$$P=P(L_1)^{\oplus m_1}\oplus\dots\oplus P(L_k)^{\oplus m_k},$$
where $m_i=\operatorname{dim}\operatorname{Hom}_{\mathfrak g}(P,L_i)$.
Similarly, any injective module $I$ in $\mathcal F(\mathfrak g)$
has a direct sum decomposition
$$I=I(L_1)^{\oplus m_1}\oplus\dots\oplus I(L_k)^{\oplus m_k},$$
where $m_i=\operatorname{dim}\operatorname{Hom}_{\mathfrak g}(L_i,I)$.

 By Lemma \ref{reciprocity} we have 
$$\operatorname{Ind}(S)=P(L_1)^{\oplus m_1}\oplus\dots\oplus P(L_k)^{\oplus m_k},$$
$$\operatorname{Coind}(S\otimes T)=I(L'_1)^{\oplus m_1}\oplus\dots\oplus I(L'_k)^{\oplus m_k},$$
with $P(L_i)\simeq I(L'_i)$. Thus, for every simple $L$ there exists
$L'$ such that  $P(L)\simeq I(L')$. Moreover, for any simple
$\mathfrak g_0$-module $S$ we have
\begin{equation}\label{eqrec1}
\operatorname{Hom}_{\mathfrak g}(\operatorname{Ind}(S),L)=\operatorname{Hom}_{\mathfrak g_0}(S,L)=
\operatorname{Hom}_{\mathfrak g}(L',(\operatorname{Coind}(S\otimes T)))=\operatorname{Hom}_{\mathfrak g_0}(L',S\otimes T)
\end{equation}

Observe first  that (\ref{eqrec1}) implies $P(\mathbb F)\simeq I(T)$, where $\mathbb F$
is the even trivial $\mathfrak g$-module. 
Therefore 
\begin{equation}\label{eqrec}
\operatorname{Hom}_\mathfrak g(P,\mathbb F)=\operatorname{Hom}_\mathfrak g(T,P)
\end{equation}
for any projective (and hence injective) module $P$.

Now we will prove that $L'=L\otimes T$ for any simple $L$. Assume the
opposite. Then there exists a simple $M$, not isomorphic to $L$
such that $P(L)\simeq I(M\otimes T)$. Note that 
$$\operatorname{Hom}_\mathfrak g(T,I(M\otimes T)\otimes
M^*)=\operatorname{Hom}_\mathfrak g(M\otimes T, I(M\otimes T))=\mathbb
F.$$
On the other hand,
$$\operatorname{Hom}_\mathfrak g(P(L)\otimes M^*,\mathbb F)=\operatorname{Hom}_\mathfrak g(P(L),M)=0.$$
We obtain a contradiction with (\ref{eqrec}).
\end{proof}

We conclude this subsection with two obvious observations. First,
every indecomposable object $M\in\mathcal F(\mathfrak g)$ admits a
central charge, i.e. every element of $Z(\mathfrak g)_0$ acts on $M$ as
a scalar operator. Hence $\mathcal F(\mathfrak g)$ decomposes into
direct sum of subcategories of $\mathfrak g$-modules with fixed
central charge. Thus, for the purpose of description of $\mathcal F(\mathfrak g)$ 
it is sufficient to consider $\mathfrak g$ with
one-dimensional center. Second, if $\mathfrak g''$ denote an
ideal in $\mathfrak g$ such that $\mathfrak g/\mathfrak g''$  is an even
abelian algebra, then every  $M''\in\mathcal F(\mathfrak g'')$ can
be lifted to a module $M\in\mathcal F(\mathfrak g)$ by decorating
$M''$ with a suitable grading. Thus, the restriction functor $\mathcal F(\mathfrak g)\to\mathcal F(\mathfrak g'')$
maps simple, injective or projective modules in $\mathcal F(\mathfrak g)$ to  simple, injective or projective modules respectively in
$\mathcal F(\mathfrak g'')$.

\subsection{Root and weight decompositions} We assume that $\mathfrak g$ is quasireductive.

Fix a Cartan subalgebra $\mathfrak h_0$ in   $\mathfrak g_0$. We call 
the centralizer $\mathfrak h$ of $\mathfrak h_0$ in $\mathfrak g$ a Cartan
subalgebra of  $\mathfrak g$. It is clear that $\mathfrak h$ is a
subalgebra and that $\mathfrak h_0\in Z(\mathfrak h)$. There is a finite subset 
$\Delta\subset\mathfrak h_0^*$ such that 
\begin{equation}\label{rootdecomp}
\mathfrak g=\mathfrak h\oplus\bigoplus_{\alpha\in\Delta}\mathfrak g_\alpha, 
\end{equation}
where
$$\mathfrak g_{\alpha}=\{x\in\mathfrak g |
[h,x]=\alpha(h)x\,\,\text{for all}\,\,h\in\mathfrak h_0\}.$$
As in the classical case $\alpha\in\Delta$ is called a {\it root} of
$\mathfrak g$ and $\mathfrak g_\alpha$ is called a {\it root space}.
By $W$ we denote the Weyl group of $\mathfrak g_0$. Then $W$ acts
naturally on $\mathfrak h_0^*$ and $\Delta$ is invariant under this action.

The following lemma follows immediately from the description of
quasireductive superalgebras given in Corollary \ref{str} and the
description of root decompositions for simple Lie superalgebras (see \cite{K}).
\begin{lemma}\label{rootspace} $\operatorname{dim}\mathfrak g_\alpha$
equals $(1|0),(1|1)$ or $(0|n)$ for some $n>0$. 
\end{lemma}

Let $Q\subset\mathfrak h_0^*$ denote the abelian subgroup generated by
$\Delta$. 

\begin{lemma}\label{rootlattice} A quasireductive $\mathfrak g$ is the Lie superalgebra
of some quasireductive algebraic supergroup if and only if the natural morphism
$Q\times_\mathbb Z\mathbb F\to\operatorname{span}\Delta$ is an isomorphism.
\end{lemma}
\begin{proof} Indeed, the adjoint action of $\mathfrak g_0$ in $\mathfrak g$
can be lifted to a representation of an algebraic group iff the above
condition holds. The lemma follows from Theorem \ref{exist}.
\end{proof} 

If $M\in\mathcal F(\mathfrak g)$, then $\mathfrak h_0$ acts semisimply
on $M$. Hence there exists a finite set $\mathcal P(M)\subset \mathfrak h_0^*$
such that
$$M=\bigoplus_{\mu\in\mathcal P(M)} M_{\mu},$$
where 
$$M_{\mu}=\{m\in M |
hm=\mu(h)m\,\,\text{for all}\,\,h\in\mathfrak h_0\}.$$
The spaces $M_\mu$ are called {\it weight spaces} and $\mathcal P(M)$ is called the
set of weights of $M$. It is clear that
$\mathcal P(M)$ is $W$-invariant.

We define the {\it character} $\operatorname{ch} M$ as a formal expression 
$$\operatorname{ch} M=\sum_{\mu\in P(M)}\operatorname{dim}M_\mu e^\mu,$$
where we use double numbers for superdimension
$$\operatorname{dim}V=\operatorname{dim}V_0+\varepsilon \operatorname{dim}V_1,$$
with $\varepsilon^2=1$. The reader can easily check that
$$\operatorname{ch}(M\oplus N)=\operatorname{ch}M+\operatorname{ch}N\,\,
\operatorname{ch}(M\otimes N)=\operatorname{ch}M\operatorname{ch}N.$$

Each weight space $M_\mu$ is an $\mathfrak h$-module. Our next
step is to describe all simple $\mathfrak h$-modules.

Let $\lambda\in\mathfrak h_0^*$. It induces a skewsymmetric bilinear
form $\omega_{\lambda}$ on $\mathfrak h_1$ defined by
$$\omega_\lambda(x,y)=\lambda([x,y]).$$
Let $\mathfrak u_1\subset\mathfrak h_1$ be a maximal isotropic subspace
with respect to  $\omega_{\lambda}$ and $\mathfrak u:=\mathfrak
h_0\oplus \mathfrak u_1$. Then $\mathfrak u$ is a subalgebra, and it is
not hard to see that $\lambda$ extended by zero on $\mathfrak u_1$ is
a character on $\mathfrak u$. Denote by $\mathbb F_\lambda$ the
corresponding even 1-dimensional representation and let
$$C_\lambda:=U(\mathfrak h)\otimes_{\mathfrak u}\mathbb F_\lambda.$$

\begin{lemma}\label{cartan} 

(1) $C_\lambda$ is a simple $\mathfrak h$-module;

(2) $C_\lambda\simeq \Pi(C_\lambda)$ iff $\operatorname {dim}
   \mathfrak h_1/\mathfrak u_1$ is odd;

(3) Every simple $M\in\mathcal F(\mathfrak h)$ is isomorphic to  
$C_\lambda$ or $\Pi(C_\lambda)$.
\end{lemma}

\begin{proof} Let $B_\lambda=U(\mathfrak h)/I_\lambda$, where
$I_\lambda$ is the ideal generated by $h-\lambda(h)$ for all
$h\in\mathfrak h_0$. 
Then $B_\lambda$ is isomorphic to the Clifford
algebra 
$T(\mathfrak h_1)/(xy+yx-\omega_\lambda(x,y))$.
The Jacobson radical of  $B_\lambda$ is generated by $\operatorname{Ker}\omega_\lambda$. 
The semisimple quotient is isomorphic to the matrix algebra of size 
$2^n$ if $\operatorname{dim}\mathfrak h_1/\operatorname{Ker}\omega_\lambda=2n$ and the direct sum of two
matrix algebras of size $2^n$ if  $\operatorname{dim}\mathfrak h_1/\operatorname{Ker}\omega_\lambda=2n+1$. Now we can
use the theory of Clifford algebras.
If  $\operatorname{dim}\mathfrak h_1/\operatorname{Ker}\omega_\lambda=2n$, 
$B_\lambda$ has one up to isomorphism simple module $C_\lambda$ (if
one disregards $\mathbb Z_2$-grading). In the category of $\mathbb Z_2$ graded modules 
$C_\lambda$ and $\Pi(C_\lambda)$ are not isomorphic since $C_\lambda^{\mathfrak u_1}$
and $\Pi(C_\lambda)^{\mathfrak u_1}$ have different parity.
If  $\operatorname{dim}\mathfrak h_1/\operatorname{Ker}\omega_\lambda=2n+1$, 
$C_\lambda$ splits into
a direct sum of two non-isomorphic simple submodules but these submodules
are not homogeneous with respect to the $\mathbb Z_2$-grading. Hence $C_\lambda$ is
simple in  the category of $\mathbb Z_2$ graded  $B_\lambda$-modules
and  $C_\lambda\simeq  \Pi(C_\lambda)$.

The above arguments imply (1) and (2). To show (3) observe that if $M$
is simple, then $\mathfrak h_0$ acts via some character $\lambda$ on $M$.
Hence $M$ is a $B_\lambda$-module, and the statement follows from the
representation theory of Clifford algebras as above. 
\end{proof}

\subsection {Highest weight theorem} Let $\gamma:Q\to\mathbb R$ be a
homomorphism of abelian groups such that $\gamma(\alpha)\neq 0$ for
any $\alpha\in\Delta$. Set
$$\Delta^\pm:=\{\alpha\in\Delta | \pm \gamma(\alpha)>0\} ,\,\,\mathfrak
n^\pm:=\bigoplus_{\alpha\in\Delta^\pm}\mathfrak g_\alpha,\,\,\mathfrak
b:=\mathfrak h\oplus \mathfrak n^+.$$
Then $\Delta=\Delta^+\cap\Delta^-$, $\mathfrak n^\pm$, $\mathfrak b$ are subalgebras
of $\mathfrak g$ and we have a {\it triangular decomposition}
$$\mathfrak g=\mathfrak n^-\oplus\mathfrak h\oplus\mathfrak n^+.$$
A subalgebra $\mathfrak b$ is called a Borel subalgebra of
$\mathfrak g$. It depends on a choice of $\gamma$ and in contrast with
the usual case not all Borel subalgebras of $\mathfrak g$ are
conjugate. Let us mention that such definition of Borel subalgebra was
first introduced in \cite{IO}.

Introduce a partial order on $\mathfrak h_0^*$ by
$$\mu\leq\lambda\,\,\text{iff}\,\,\lambda=\mu+\sum_{\alpha\in\Delta^+}m_\alpha\alpha,\,m_\alpha\in\mathbb
Z_{\geq 0}.$$

\begin{theorem}\label{highestweight} Let $L\in\mathcal F(\mathfrak g)$
be simple. There exists a unique weight $\lambda\in \mathcal P(L)$, called the
highest weight of $L$, such that $\mu\leq\lambda$
for all $\mu\in \mathcal P(L)$. Furthermore, $L_\lambda$ is a simple $\mathfrak h$-module.
If two simple modules $L$ and $M$ have the same highest weight, then
$M$ is isomorphic to $L$ or $\Pi(L)$.
\end{theorem}
\begin{proof} This can be proven in the same way as the analogous
result for reductive Lie algebras. So we will be brief. Pick up a
maximal element $\lambda$ in $\mathcal P(L)$.  Since $\mathfrak g_\alpha L_\mu\subset L_{\mu+\alpha}$,
$\mathfrak n^+ L_\lambda=0$. Therefore by PBW theorem
$L=U(\mathfrak n^-)L_\lambda$, and  
$\mu\leq\lambda$ for every $\mu\in\mathcal P(L)$. Every proper $\mathfrak h$-submodule of $L_\lambda$ 
generates a proper submodule of $L$. Thus, the simplicity of $L$
implies the simplicity of $L_\lambda$. 

Finally,  Lemma \ref{cartan} implies that
$L_\lambda\simeq C_\lambda$ or $\Pi(C_\lambda)$. Define a Verma module $M(\lambda)=U(\mathfrak g)\otimes_{U(\mathfrak b)}C_\lambda$. 
It is easy to see that $M(\lambda)$ has a unique proper maximal
submodule and hence a unique simple quotient. So $L$ is either the
simple quotient of $M(\lambda)$ or $\Pi(M(\lambda))$. That implies the
last statement of the theorem.
\end{proof}
\begin{corollary}\label{char} If $L,M\in\mathcal F(\mathfrak g)$ are
simple and $\mathcal P(L)=\mathcal P(M)$, then $M\simeq L$ or $\Pi(L)$.
\end{corollary}

It was first shown in \cite{K2} that one can write a superanalogue of
the Weyl formula for $\operatorname{ch} L$ if $L$ has a generic
highest weight. This result was pushed further in \cite{Prep}. In
particular, it can be done for all quasireductive superalgebras. On
the other hand, the problem of finding $\operatorname{ch} L$ for all
finite-dimensional simple $L$ is solved only
for $\mathfrak g=\mathfrak{gl}(m,n),\mathfrak{osp}(m,2n)$ and
$\mathfrak q(n)$,
see  \cite{S},\cite{CLW},\cite{GrS},\cite{B2},\cite{PS1}.

A reasonable description of the category $\mathcal F(\mathfrak g)$ is rather complicated task. The case of general linear
and orthosymplectic superalgebras is done in \cite{BS} and \cite{CLW}.
In general, this is an open problem.

\begin{example}
Let $\mathfrak g=\tilde{\mathfrak k}^d$. It is not difficult to see that any non-trivial simple module
$M\in\mathcal F(\mathfrak g)$ is isomorphic to $M_0\otimes \mathbb F(\theta)$, where $M_0$ is a simple 
$\mathfrak k$-module (maybe with shifted action $\theta\frac{\partial}{\partial \theta}$). 
The projective and injective indecomposable modules have the following interpretation. Consider
the chain complex $\mathcal C(\mathfrak k,M_0)=\Lambda(\mathfrak k)\otimes M_0$ with the
differential $d$, calculating the homology $H_*(\mathfrak k,M_0)$. 
For any $x\in\mathfrak k, c\in\mathcal C(\mathfrak k,M_0)$ let 
$m_x(c)=x\wedge c$. Then $d$, $m_x$ and the grading  operator generate
the Lie superalgebra isomorphic to  $\mathfrak g$. We leave to the
reader to check that $\mathcal C(\mathfrak k,M_0)$ is isomorphic to
the projective cover $P(M)$ in $\mathcal F(\mathfrak g)$.
Similarly $I(M)$ can be identified with the cochain complex calculating
the cohomology $H^*(\mathfrak k, M_0)$.
\end{example}

\begin{example} Let $\mathfrak g$ be a quasireductive pseudoabelian
superalgebra with even 1-dimensional center $\mathbb Fz$. Every
indecomposable $\mathfrak g$-module $M$ admits a central charge $c$.
If $c=0$ then every simple module $M$ is simple over $\mathfrak g_0$
with trivial action of $\mathfrak g_1$, and
$P(M)=\operatorname{Ind}(M)$, $I(M)=\operatorname{Coind}(M)$.

Let $c\neq 0$. Consider the skewsymmetric form $(\cdot,\cdot)$ on $\mathfrak g_1$
defined by $[x,y]=c (x,y)$. For simplicity assume that $(\cdot,\cdot)$
is non-degenerate. Let $\operatorname{Cliff}(\mathfrak g_1)$ be the
Clifford algebra defined by the form  $(\cdot,\cdot)$. Then we have an
embedding $\mathfrak g\to\operatorname{Cliff}(\mathfrak g_1)$.
A simple $\operatorname{Cliff}(\mathfrak g_1)$-module $V$ is
simple as a $\mathfrak g$-module. Every simple $\mathfrak g$-module $M$ with
central charge $c\neq 0$ is isomorphic to $M'\otimes V$ for some
simple $\mathfrak g$-module $M'$ with zero central charge. It is easy to see that $M$ is projective
and injective, and hence the subcategory of $\mathfrak g$-modules with
non-trivial central charge is semisimple. 
\end{example}

\begin{example} Let $\mathfrak g=p(n)$ with $n>2$. Consider
the $\mathbb Z$-grading 
$$\mathfrak g^{-1}\oplus\mathfrak g_0\oplus\mathfrak g^1,$$
where $\mathfrak g_0=\mathfrak {gl}(n)$ and $\mathfrak g^{1}=S^2(E)$,
 $\mathfrak g^{-1}=\Lambda^2(E^*)$, $E$ being the standard
 $\mathfrak {gl}(n)$-module. Let 
$$\mathfrak g^{\pm}:=\mathfrak g_0\oplus\mathfrak g^{\pm 1}$$ and
$$K^\pm(\lambda):=U(\mathfrak g)\otimes_{U(\mathfrak g^\pm)}L^0(\lambda),$$
where $L^0(\lambda)$ is the simple $\mathfrak g_0$-module with highest
weight $\lambda$ and trivial
$\mathfrak g^{\pm 1}$-action. It is not hard to see $K^\pm(\lambda)$
is indecomposable with unique simple submodule and unique simple quotient.
Note that $K^-(\lambda)$ is never simple, but $K^+(\lambda)$ is simple
for generic $\lambda$. It was shown in \cite{K2} that, if $\lambda=(a_1,\dots,a_n)$ 
in the standard $SL(n)$ notation,  $K^+(\lambda)$ is simple iff 
$\prod_{i\leq j}(a_i-a_j)\neq 0$. A simple module is never projective.
A description of the socle filtration in $K^\pm(\lambda)$ and 
calculation of the character of a simple $p(n)$-module are open problems.

If $n=4$,  $\mathfrak g^{-1}=\Lambda^2(E^*)$ admits an $\mathfrak{sl}(4)$-invariant
symmetric form. That implies the isomorphism  $\mathfrak{sl}(4)\simeq\mathfrak{so}(6)$. By this reason $sp(4)$ has a non-trivial central
extension $\hat{sp}(4)$. The representation theory of  $\hat{sp}(4)$
with non-zero central charge is quite different from that in the case of $p(4)$.
Although $K^+(\lambda)$ can be defined in this case and has the same dimension,
the condition of its simplicity should be different. Note that
$\mathfrak f=\hat{sp}(4)$ has a $\mathbb Z$-grading
$$\mathfrak f=\mathfrak f^{-2}\oplus \mathfrak f^{-1}\oplus\mathfrak f^{0}\oplus\mathfrak f^{1},$$
with $\mathfrak f^0\simeq \mathfrak{so}(6)$,  $\mathfrak f^{-1}$ being
the standard representation of $\mathfrak{so}(6)$. The subalgebra
$\mathfrak f^-:=\mathfrak f^{-2}\oplus \mathfrak f^{-1}\oplus\mathfrak f^{0}$ 
is pseudoabelian and hence  we have an embedding  $\mathfrak f^-\to\operatorname{Cliff}(\mathfrak f^{-1})$.
This embedding can be extended to $\mathfrak f\to\operatorname{Cliff}(\mathfrak f^{-1})$. 
The restriction of a simple $\operatorname{Cliff}(\mathfrak f^{-1})$
to $\mathfrak f$ is an example of a simple $(4|4)$-dimensional $\mathfrak f$-module with non-trivial central charge. 
\end{example}

\end{document}